\documentclass[reqno,centertags]{amsart}
\usepackage{amsmath,amsthm,amscd,amssymb}
\usepackage{latexsym}
\usepackage{bbm}

\newcommand{\lb}{\label}

\newcommand{\dom}{\text{\rm{dom}}}

\newcommand{\beq}{\begin{equation}}
\newcommand{\eeq}{\end{equation}}
\newcommand{\ba}{\begin{align}}
\newcommand{\ea}{\end{align}}

\def \d{{\tt d}}

\def \sur#1#2{\mathrel{\mathop{\kern 0pt#1}\limits^{#2}}}
\def \el{\sur{=}{(d)}}

\newcounter{smalllist}

\newcommand{\Pois}{\operatorname{Pois}}

\allowdisplaybreaks
\numberwithin{equation}{section}

\newtheorem{theorem}{Theorem}[section]
\newtheorem{proposition}[theorem]{Proposition}

\theoremstyle{definition}

\theoremstyle{remark}
\newtheorem{remark}{Remark}

\begin{document}
\title[Invariance principles for clocks]
{Invariance principles for clocks}
\author[M-E.~Caballero and A.~Rouault]
{Maria-Emilia Caballero$^{1}$ and
Alain Rouault$^{2}$}

\thanks{$^1$ Instituto de Matem\'aticas, UNAM, Coyac\'acan  04510, M\'exico DF, Mexique, email: marie@matem.unam.mx}

\thanks{$^2$ Laboratoire de Math{\' e}matiques de Versailles, UVSQ, CNRS, Universit\'e Paris-Saclay, 78035-Versailles Cedex France, e-mail: alain.rouault@uvsq.fr}

\date{\today}

\begin{abstract} 
We show an invariance principle for rescaled clocks of positive semi-stable Markov processes, proving a conjecture presented in Remark 4 in  Demni, Rouault, Zani \cite{demni2015large}, 2015.
\end{abstract}

\maketitle

\section{Introduction - pssMp and OU} \lb{s1}

For $\alpha > 0$, a 
 positive self-similar 
 Markov process (pssMp) 
  of index $\alpha$,   is a $[0, \infty)$-valued strong Markov process $(X, \mathbb Q_a), a > 0$ with c\`adl\`ag paths, fulfilling the scaling property
\begin{equation}
\label{self}
\left(\{b X_{b^{-\alpha}t}, t \geq 0\}, \mathbb Q_a\right) \el \left(\{X_t , t \geq 0\}, \mathbb Q_{ba}\right) 
\end{equation}
for every $a, b > 0$.

The Lamperti transformation (see \cite{lamp}) connects these processes to L\'evy processes.  Let us summarize this connection. We will follow the notations of \cite{demni2015large}.

Any pssMp $X$   which never reaches the boundary state $0$
 may be expressed as the
exponential of a L\'evy process  not drifting to $-\infty$, time changed by the inverse of its exponential functional. 
More formally, if $(X, (\mathbb Q_a)_{a > 0})$ is  a  pssMp of index $\alpha$  which never reaches $0$,  set
\begin{equation}
T^{(X)}(t) = \int_0^t \frac{ds}{X_s^\alpha} \  , \ (t \geq 0)
\end{equation}
and let $A^{(X)}$ be  its inverse, defined by
\begin{equation}
A^{(X)}(t) = \inf \{u \geq 0 : T^{(X)}(u) \geq t \} \  , \ (t \geq 0)\,,
\end{equation}
and let $\xi$ be the process defined by
\begin{equation}
\xi_t = \log X_{A^{(X)}(t)} - \log X_0\ , \ (t\geq 0)\,.
\end{equation}
Then, for every $a >0$, the distribution of $(\xi_t, t \geq 0)$ under $\mathbb Q_a$ does not depend on $a$ and is the distribution of a L\'evy process starting from $0$.

Conversely, let  $(\xi_t, t \geq 0)$ be a  L\'evy process starting from $0$ and let $\mathbb P$ and $\mathbb E$ denote the underlying probability and expectation, respectively.  

 Fix $\alpha > 0$. Set
\begin{equation}
\label{expfunc}
\mathcal{A}^{(\xi)} (t) = \int_0^t e^{\alpha \xi_s}ds\,.
\end{equation}
Here, we  assume that  $\xi$ does not drift to $-\infty$ i.e. satisfies,   
$\limsup_{t\uparrow \infty} \xi_t = \infty$. The
 inverse process $\tau^{(\xi)}$ of $\mathcal A^{(\xi)}$ is
\begin{equation}
\label{invexpfunc}
\tau^{(\xi)}(t) = \inf\{u \geq 0 : \mathcal{A}^{(\xi)}(u) \geq t\} \  , \ (t \geq 0)\,.
\end{equation}

For every $a >0$, let $\mathbb Q_a$ be the law under $\mathbb P$ of the time-changed process

\begin{equation}
\label{xitoX1}
X_t = a\exp \xi_{\tau^{(\xi)}(ta^{-\alpha})} \ , \ (t\geq 0)\,,
\end{equation}
then $(X , (\mathbb Q_a)_{a > 0})$ is a  pssMp of index $\alpha$  which never reaches $0$ and we have the fundamental relation
\begin{equation}
\label{fund}
\tau^{(\xi)} (t) = T^{(X)} (t X_0^\alpha)\,.
\end{equation}

The process $(T^{(X)} (t) , t\geq 0)$ is called the \underline{clock} associated with the pssMp $X$.
Some years ago, a particular interest was dedicated to the asymptotic behaviour of this process in long time (\cite{zani}, \cite{demni2015large}). To recall the Law of Large Numbers we need some notations.

Let $\psi$ be the Laplace exponent of $\xi$, defined by
\begin{align}
\label{defexp}\mathbb E \exp (m\xi_t) = \exp (t\psi(m))\,,
\end{align}
and set $\dom\  \psi = \{ m : \psi(m) < \infty\}$.
We assume
\begin{align}
\label{assume}
0\in \hbox{int}\ \dom\ \psi \  \hbox{and} \  
 p := \mathbb E \xi_1 = \psi'(0) > 0\,.
\end{align}

Starting from $a >0$ such a process
 $X_t$ never hits $0$. Nevertheless, a probability measure $\mathbb Q_0$
 can be obtained as the weak limit of $\mathbb Q_a$ when $a\downarrow 0$ and under $\mathbb Q_0$  the canonical
process has the same transition semigroup as the one associated with $(X; (\mathbb Q_a)_{a>0}$. 
(see \cite{bertoin2002entrance}, \cite{bertcaba}, \cite{caballero2006weak}, \cite{chaumont2012fluctuation}, \cite{pardo2013self}). A sufficient condition is (\ref{assume}) plus
\begin{align}
\label{arith}
\hbox{the support of $\xi$ is not arithmetic}\,.
\end{align} 
This latter measure is an entrance law for the
 semigroup $p_t f(x) = \mathbb E_x f(X_t)$ 
and satisfies
\begin{equation}
\label{magicI}
\mathbb Q_0(f(X_t)) = \frac{1}{\alpha\mathbb E \xi_1} \mathbb E\left[I_\infty^{-1}f\left((t/I_\infty)^{1/\alpha}\right)\right]
\end{equation}
where 
\[I_\infty 
 = \int_0^\infty e^{-\alpha\xi_s}  ds\,.\]
 The Law of Large Numbers is the following.

\begin{theorem}[\cite{demni2015large} Th.1]
\label{thmLLN}
Assume (\ref{assume}) and (\ref{arith}). As $t \rightarrow \infty$,
\begin{enumerate}
\item For every $a >0$, 
\begin{equation}
\label{LLNX}
\frac{1}{\log t} \int_0^t \frac{ds}{X_s^\alpha}\rightarrow (\alpha p)^{-1} 
 \,, \ \ \mathbb Q_a-\hbox{a.s.}
\end{equation}
\item
\begin{equation}
\label{LLN0}
\frac{1}{\log t} \int_1^t \frac{ds}{X_s^\alpha}\rightarrow (\alpha p)^{-1} 
 \,, \ \ \mathbb Q_0-\hbox{a.s.}
\end{equation}
\end{enumerate}
\end{theorem}
In \cite{demni2015large}, the authors go on with the study of large deviations. 
A CLT for Bessel clocks was previously  proved in Exercise X.3.20 in \cite{RY} and extended to an invariance principle in (\cite{zani}), using stochastic analysis.  Following the result on LDP for clocks, a general CLT is conjectured in \cite{demni2015large} Remark 4. In Section 2, we state an invariance principle (Functional CLT) for this kind of processes.
proved in  Section 3 and illustrated by examples in Section 4.

The main tool is the introduction of an  ergodic process with nice asymptotic properties. 
With a pssMp  $(X_t)$ of index $\alpha$, it is classical to associate a process called generalized Ornstein-Uhlenbeck (OU) by
\begin{align}
\label{OU1def}
U(t) = e^{-t/\alpha} X (e^t)
\end{align}
which 
is strictly stationary, Markovian, ergodic under $\mathbb Q_0$, and its invariant measure is the law of $X_1$ under $\mathbb Q_0$ i.e.
\begin{align}
\label{invlaw}
\mu^U(f) =\frac{1}{\alpha\mathbb E\xi_1}\mathbb E\left[I_\infty^{-1} 
f\left( I_\infty^{-1/\alpha}\right)\right]\,.
\end{align} 
The infinitesimal generators $L^X$ and $L^\xi$ are related by
\begin{align}
\label{genX}
L^X h(x) = x^{-\alpha}L^\xi (h \circ \exp)(\log x)
\end{align}
(see \cite{CPY2}, \cite{pardo2013self}), and the generator of $U$ is 
\begin{align}
\label{genU}
L^U h (x) =  L^X h(x) -\frac{x}{\alpha}h'(x)\,.
\end{align}
Two examples of function $f$ are simple and particularly useful.

If $m \in \dom\ \psi$   the function  $x \mapsto \exp (m x)$ is in the domain of $L^\xi$ , so that 
\[f_m : x \mapsto x^m\,,\]   is in the domain of $L^X$ and we have 
\begin{align}L^X f_m (x) = \psi (m) f_{m-\alpha}\,, \end{align}
so that
\begin{align}
\label{fmOU}L^U f_m =  \psi (m) f_{m-\alpha} - \alpha mf_m\,.\end{align}
The  function $\mathfrak i (x) = x$ is in the domain of $L^\xi$ with $L^\xi \mathfrak i = p$.  Formula  (\ref{genX}) tells us that the function $\varphi(x) = \log x$ is in the domain of $L^X$ and
\begin{align}
L^X \varphi (x) = \frac{1}{x^\alpha}L^\xi \mathfrak i\!\ (\log x) 
= \frac{p}{x}\,, 
\end{align}
and owing to (\ref{genU})
\begin{align}
\label{logOU}
L^U \varphi (x) = \frac{p}{x^\alpha} -\frac{1}{\alpha}\,.
\end{align}
\begin{remark}
\label{2OU}
 There is a variant of $U$, defined by
\[\tilde U(t) = e^{-t/\alpha}X(e^t -1)\]
which shares the same transition with $U$ and begins at $X(0)$ at $t=0$ (\cite{bertoin2019ergodic}).
\end{remark}
\begin{remark}
\label{scale}
If $(X, (\mathbb Q_a)_{a >0})$ is  a pssMa of index $\alpha$, then  the process $Y =( X^{\alpha}, (\mathbb Q_{a^\alpha})_{a >0})$, is a pssMp of index $1$.
 Conversely if $(Y, (\mathbb Q_a)_{a > 0})$ is a pssMp if index $1$
then, for any $\alpha > 0$, 
the process 
$(X = Y^{1 / \alpha} ,(\mathbb Q_{a^{1/\alpha}})_{a > 0})$
 is a pssMp of index $\alpha$.
\end{remark}
\section{Main result}
The following theorem states an invariance principle  ( or functional central limit theorem FCLT) under two regimes,  $\mathbb Q_0$ and (under conditions) $\mathbb Q_a, a >0$.
\begin{theorem}
\label{main}
\begin{enumerate}
\item
Under $\mathbb Q_0$,  as $T \rightarrow \infty$,
\begin{align}
\label{underQ0}
\left( (\log T)^{-1/2}\left(\int_1^{T^t} \frac{dr}{X^\alpha(r)} - \frac{t\log T}{\alpha p}\right) ; t \geq 0\right)\Rightarrow \left(v W(t) ;  t \geq 0\right)\,.
\end{align}
where
\begin{align}
v^2 =\frac{\sigma^2}{\alpha p^3}, \ \sigma^2 
 = \psi"(0)\,. 
\end{align}
\item
If one of the following conditions
\begin{enumerate}
\item
there exists $m>0$ such that $\psi(m) <0$ 
\item there exists $m > \alpha$ such that $\psi(m) > 0$,
\end{enumerate} is satisfied,  then for every $a >0$, under $\mathbb Q_a$, as $T \rightarrow \infty$, 
(\ref{underQ0}) holds true.
\end{enumerate}
\end{theorem}
\begin{remark}
When the above criterion is not checked, the invariance principle holds true under $\mathbb Q_a$ for almost every $a$ (see Theorem 2.8 in \cite{bhattacharya1982functional} ).
\end{remark}

\section{Proof of the main result}
\subsection{FCLT under the invariant measure}
Observe that, owing to (\ref{OU1def})
\begin{align}
\label{convXU}
\int_1^{T^t} \frac{ds}{(X(s))^\alpha} - \frac{t \log T}{\alpha p} = \int_0^{t\log T} \left(\frac{1}{U(r)^\alpha} - \frac{1}{\alpha p}\right) dr
\end{align}
which reduces the problem to an invariance principle for a functional of the process $U$.
We will use a classical result on weak convergence.
\begin{theorem}[Bhattacharya  Th. 2.1 \cite{bhattacharya1982functional}]
\lb{B1}
Let $(Y_t)$ be a measurable stationary ergodic process with an invariant probability $\pi$. If $f$ is in the range $\hat A$ of the extended generator 
 of $(Y_t)$, then, as $n \to \infty$ 
\begin{align}
\lb{C1}\left(n^{-1/2}\int_0^{nt} f(Y_s) ds\right)_{t \geq 0} \Rightarrow (\rho W(t))_{t \geq 0}
\end{align}
where
\begin{align}
\label{-2}\rho^2 = - 2 \int f(x) g(x) \pi(dx) \ , \  \hat A g= f\,.
\end{align}
\end{theorem}

\noindent Owing to (\ref{logOU}) we see that the pair $(f,g)$ with
\begin{align}
\label{pair}
f(x) = \frac{1}{x^{\alpha}} - \frac{1}{\alpha p} \ , \ g(x) = \frac{1}{p} \log x\,,\end{align}
satisfies
\begin{align}
\label{LUg}
L^U g = f\,.
\end{align}
Now the convergence in distribution comes from  Th. \ref{B1}  and formula (\ref{convXU}), with
\begin{align}
\label{defv}
v^2 = - 2 \int f(x) g(x) \mu^U(dx)
\end{align}
where $\mu^U$ is the invariant distribution of the process $U$.

It remains to compute the variance $v^2$. From (\ref{invlaw}) and (\ref{defv}) we have
\begin{align}
\notag
v^2 &=- 2 (\alpha p)^{-1}\mathbb E \left(I_\infty^{-1} f(I_\infty^{-1}) g(I_\infty^{-1})\right)\\\notag &= 2(\alpha p)^{-1} \mathbb E \left(I_\infty^{-1}\left(I_\infty - (\alpha p)^{-1}\right)p^{-1} \log  I_\infty\right)\\&= 2(\alpha p)^{-2}\mathbb E \left(\log I_\infty -  (\alpha p)^{-1}I_\infty^{-1}\log I_\infty\right)\,.
\end{align}
The Mellin transform 
\[M(z) = \mathbb E(I_\infty^{-z})\]
may play a prominent role, since
\[\mathbb E(\log I_\infty^{-1}) =  M'(0) \ ; \  \mathbb E \left(I_\infty^{-1}\log I_\infty^{-1}\right) =  M'(1)\,,\]
so that 
\begin{align}
v^2 = 2 (\alpha p)^{-2}\left(-M'(0) + (\alpha p)^{-1}M'(1)\right)\,.
\end{align}
We only know that $M$ satisfies the recurrence equation
\begin{align}
\label{recursion}
\psi(\alpha z) M(z) = z M(z+1)\,.
\end{align}
(see \cite{bertoin2004exponential} Th. 2 i) and Th. 3. and apply scaling).
Differentiating twice the above formula gives
\[ \alpha^2\psi''(\alpha z)M(z) +2 \alpha \psi'(\alpha z) M'(z) + \psi(\alpha z) M''(z) = 2 M'(z+1) + z M''(z+1)\]
which, for $z=0$, gives
\[\alpha^2\psi''(0) + 2 \alpha p  M'(0) = 2 M'(1)\]
and then
\begin{align}
v^2 = \frac{\psi''(0)}{\alpha p^3} = \frac{\sigma^2}{\alpha p^3}\,.
\end{align}
This fits exactly with  the conjecture in \cite{demni2015large} Rem. 4.
\subsection{Quenched FCLT}
We want to show 
the  FCLT  under $\mathbb Q_a$ for every $a > 0$.

\begin{theorem}[Bhattacharya Th. 2.6 \cite{bhattacharya1982functional}]
Let $(p_t)_{t \geq 0}$ be the semigroup of a Markov process $(Y_t)$. Assume that for every $x$, as $t \to \infty$
\[\Vert p_t (x ; \cdot ) - \mu\Vert_{\hbox{var}} \to 0\,.\]
Then, with the notations of Theorem \ref{B1},
the convergence (\ref{C1}) holds under $P_x$ for every $x$. 
\end{theorem}

Such a process  is called positive recurrent. Moreover, the exponential ergodicity is defined by the existence of a finite function $h$ and a constant $\gamma$ such that for every $x$
\[\Vert p_t(x, .) -\pi \Vert_{var} \leq h(x) e^{-\gamma t}\,.\]

Let us first examine the possibility of a general criterion on 
 the exponent $\psi$ such that the assumptions of the latter theorem are fulfilled.
A sufficient condition is given by the so-called Forster-Lyapounov drift criterion, due to \cite{meyn1993stability}).
\begin{theorem}[Wang Th. 2.1  \cite{wang2008criteria}]
\label{W2}
Let  $L$ the generator of a Markov process $Y_t$, 
\begin{enumerate}
\item If there exists a continuous function satisfying
\begin{equation}
\lb{bonne}
\lim_{|x|\to \infty} f(x) = \infty\end{equation}
and constants $K > 0, C> 0, D \in (-\infty, \infty)$ such that
\begin{equation}
\label{pr}
Lf \leq - C + D1_{[-K, K]}
\end{equation}
then the process $(Y_t)$ is positively recurrent.
\item
If there exists a continuous function satisfying (\ref{bonne}) and constants $K > 0, C> 0, D \in (-\infty, \infty)$ such that
\begin{equation}
\label{ee}
Lf \leq - C f + D1_{[-K, K]}
\end{equation}
then the process $(Y_t)$ is exponentially ergodic.
\end{enumerate}
\end{theorem}

We want to  check these criteria for our models, using the function  $f_m$ and (\ref{fmOU}). 
\begin{enumerate}
\item
When $0 < m < \alpha$ and  $\psi(m) > 0$, $f_{m-\alpha}$  is not bounded in the neighbouring of $0$, and 
then is not convenient.

\item 
Let us look for $K, C, D$ in the other cases. For  every $C \in (0,\alpha m)$, let us define
\[h_m = L^U f_m + Cf_m = \psi (m) f_{m-\alpha} + (C- \alpha m)f_m \,.\]
\begin{enumerate}
\item
 If $\psi(m) < 0$, then $h_m \leq 0$ and then for $D=0$, (\ref{ee}) holds true for every $K$.

\item
 If $\psi(m) > 0$ and $m >\alpha$, then $h_m$ is increasing for 
\[0 < x <x_{\max}= \left(\frac{(m-\alpha)\psi(m)}{m(m-\alpha C)}\right)^{1/\alpha}\] and decreasing after. It is $0$ for $x= x_0= \left(\psi(m)/(m- \alpha C)\right)^{1/\alpha} > x_m$. Then choose  $K= \psi(m) / (m-C)$ and 
$D = h_m (x_m)$
and (\ref{ee}) holds true.
\end{enumerate}
\end{enumerate}
\section{Examples}
In this section, we consider examples taken from  \cite{demni2015large}. 

 When the function $\psi$ is rational, the  distribution of $I_\infty$ may be found in  \cite{kuznetsov2012distribution}.
Let us notice that in all these examples we compute the variance using the elementary formula
\[\psi''(0) = 2 \frac{d}{dm}\frac{\psi(m)}{m}\big|_{m=0}\,.\]
\subsection{Brownian motion with drift} This is also the Cox Ingersoll Ross model.
Let us consider the L\'evy process
\[\xi_t = 2B_t + 2\nu t\,,\]
where $B_t$ is the standard linear Brownian motion and $\nu > 0$. 
In this case, $X_t$ is the squared  Bessel process of dimension $d = 2(1 +\nu)$. Its index is $1$ and it is the only continuous pssMp (of index $1$).
We have
\[\psi(m) = 2m(m+\nu), \ p = 2\nu, \ \sigma^2 = 4\ , \ v^2=(4\nu^3)^{-1} ,.\]
Condition (\ref{ee}) of Th. \ref{W2} is satisfied (see(2) b) above),  we have exponential ergodicity, the FCLT holds under $\mathbb Q_a$ for every $a\geq 0$.

The invariant measure for $U$ is easy to determinate, since 
\[I_\infty^{-1} \el 2\mathcal Z_\nu\]
where $\mathcal Z_\nu$ is gamma with parameter $\nu$ 
(see \cite{bertoin2002entrance} (6)), so that
\[\mathbb Q_0(\varphi) = \frac{1}{\mathbb E \xi_1}\int_0^\infty y 
\varphi(y)  \frac{y^{\nu-1}}{2^\nu} e^{-y/2} dy\]
which means that the invariant measure is the distribution of $2\mathcal Z_{\nu +1}$. 

\begin{remark}
Actually
it is proved in \cite{li2019exponential} Rem. 1.2 (2) that  $U$   is exponentially ergodic but not strongly ergodic, where strong ergodicity is defined 
by the existence of $\gamma > 0$ such that
\[ \sup_{x}  \Vert P_t(x, .) -\pi \Vert_{var} \leq e^{-\gamma t}\,.\]
\end{remark}

\subsection{Poissonian examples}
Let $\Pois(a,b)_t$ be the compound Poisson process of parameter $a$ whose jumps are exponential of parameter $b$.
We will consider three models : $\xi_t = \d t +\Pois(a,b)_t$ with $\d >0$, $\xi_t= - t + \Pois(a,b)_t$ and $\xi_t = t-\Pois(a,b)_t$. 

\subsubsection{$\xi_t = \d t + \Pois(a,b)_t$}
\lb{A}
\[\psi(m) = m\left(\d + \frac{a}{b -m}\right) \ \ (m\in (-\infty,  b))\,,\]
\[p = \d + \frac{a}{b}\ , \ \sigma^2 = \frac{a}{b^2}\ , \ v^2 = \frac{ab^3}{a + \d b)^3}\,.\]

\subsubsection{$\xi_t = - t + \Pois(a,b)_t$ with $b <a$}
\lb{B}
\[\psi(m) =  m\left(-1 +  \frac{a}{b -m}\right)  \ \ (m\leq b)\,,\]
\[p = \frac{a-b}{b}\ , \ \sigma^2 = \frac{a}{b^2} \ , \ v^2 = \frac{ab}{(a-b)^3}\,.\]

When $\delta = \d >0$, $I_\infty \el \alpha^{-1}B(1+ b, a \alpha^{-1})$  where $B(u,v)$ is the Beta distribution  of parameters $(u,v)$ (see \cite{gjessing1997present} Th. 2.1 i)).
The invariant measure is then the distribution of $W^{-1}$ where $W \el \alpha^{-1} B(b, a \alpha^{-1})$.

When $\delta= -1$, $I_\infty \el B_2(1+b , a-b)$ (see \cite{gjessing1997present} Th. 2.1  j)), where $B_2(u,v)$ is the Beta distribution of the second order of parameters $(u,v)$. The invariant measure is then  the $B_2 (a-b+1, b)$ distribution.
\subsubsection{$\xi_t = t  -\Pois(a,b_t)$ with $b > a$}
\lb{C} 
It is the so called spectrally negative saw-tooth process.
\[\psi(m) = m\left(1 - \frac{a}{b + m}\right) \ , \ (m\in (- b, \infty))\,,\]
\[p = \frac{b-a}{b}\ , \ \sigma^2 = \frac{a}{b^2}\ , \ v^2 = \frac{ab}{(b-a)^3}\,.\]
 We have $I_\infty^{-1}\el B(b-a, a)$  (see \cite{kuznetsov2012distribution} Th. 1), so that the invariant measure is
$B(b-a + 1, a)$.
\medskip

For examples \ref{A} and \ref{B}, $f_m$ is in the domain of the generator iff $m < b$; for example \ref{C} , $f_m$ is always in the domain. 

Looking at our above criterion, we see that we have the exponential ergodicity in the first two cases when $b >1$, and completely for the third one.
\subsection{Spectrally negative process conditioned to stay positive}
For $\alpha \in (1,2)$, let $X^{\uparrow}$ be the spectrally $\alpha$-stable process conditioned to stay positive as defined in 
 \cite{caballero2006conditioned} Sect. 3.2. and \cite{patie1} Sect. 3.
Its corresponding L\'evy process has Laplace exponent
\[\psi(m) =  \frac{\Gamma(m+\alpha)}{\Gamma(m)}\ , \  (m \in (-\alpha, \infty))\,.\]
We have
\[p = \Gamma(\alpha)\ , \  \sigma^2=
 2\left( \Gamma'(\alpha) +\gamma  \Gamma(\alpha)\right)\]
where $\gamma = - \Gamma'(1)$ is the Euler constant,
and then
\[v^2 = 2 \frac { \Gamma'(\alpha) + \gamma\Gamma(\alpha)}{\alpha \Gamma(\alpha)^3}\,.\]
Using (\ref{recursion}) one sees directly that
$M(z) = \Gamma(\alpha z+1)/\Gamma(z+1)$, so that
$I_\infty \el S_{1/ \alpha}(1)$, the stable subordinator of index $1/\alpha$ evaluated at $1$. 

Condition (2) b) of Th. 2.1 is satisfied.

\subsection{Hypergeometric stable process}

The modulus of a Cauchy process in $\mathbb R^d$ for $d > 1$ is a 1-pssMp  with infinite lifetime. The associated L\' evy process is a particular case of hypergeometric stable process of index $\alpha$ as defined in \cite{cab}, with $\alpha < d$. The characteristic exponent given therein by  Th. 7 yields the Laplace exponent :
\[\psi(m) = - 2^\alpha \frac{\Gamma((- m +\alpha)/2)}{\Gamma(-m/2)} \frac{\Gamma((m+d)/2)}{\Gamma((m+d-\alpha)/2)}\ , \ (m\in (-d, \alpha))\,.\]
We have
\[p = 2^{\alpha-1}\frac{\Gamma(\alpha/2) \Gamma(d/2)}{\Gamma(d-\alpha)/2)} \ , \ \sigma^2 = p \left[1 - \gamma - \Psi((d-\alpha)/2) - \Psi (\alpha/2)\right]\,,\]
where $\Psi$ is the Digamma function.

The distribution of the limiting variable 
$I_\infty$ is studied in \cite{kuznetsov2013fluctuations} and \cite{kuznetsov2012distribution}. 

Condition (2) b) of Th. 2.1 is never satisfied. Condition (2) a) can be satisfied if $\alpha > 2$, taking $m \in (2, \min(\alpha, 4))$ since
$\Gamma(-m/2) > 0$ hence $\psi(m) < 0$.

\subsection{Continuous State branching process with immigration (CBI)}

Let $\kappa \in [0,1)$ and $\delta> \kappa/(\kappa +1)$. Let $X$ be the continuous state branching process with immigration (\cite{Kyp} Sec. 13.5) whose branching mechanism is 
\[\phi(\lambda) = \frac{1}{\kappa}\lambda^{\kappa +1}\]
and immigration mechanism is
\[\chi(\lambda) = \delta \phi'(\lambda)\,.\]
We have the representation
\[\phi(\lambda) = \int_0^\infty(e^{-\lambda z} - 1+ \lambda z)\mu(dz) \ , \mu(dz) = \frac{\kappa +1}{\Gamma(1-\kappa)} \frac{dz}{z^{\kappa +2}}\,.\]
This process is self-similar of index $\kappa$ (see \cite{patie1}, lemma 4.8)\footnote{Beware, our $\phi$ is $-\varphi$ therein} and the corresponding Laplace exponent is
\[\psi (m) = c (\kappa - (\kappa +1) \delta -m)\frac{\Gamma(-m+\kappa)}{\Gamma(-m)}  \ , \ (m \in (-\infty, \kappa))
\,\]
and
\begin{align}
\notag
p= c ((\kappa+1)\delta- \kappa) \Gamma(\kappa) > 0 \ &, \ \sigma^2 = c\left(\Gamma(\kappa) + (\kappa - (\kappa +1) \delta) (\Gamma'(\kappa) +\gamma\Gamma(\kappa)\right)
\\ v^2 &= \frac{\sigma^2}
{\kappa p^3}\,.
\end{align}
We can then apply Th. \ref{main}(1) and conclude
that under $\mathbb Q_0$, as $T\rightarrow \infty$
\begin{align}
\left((\log T)^{-1/2} \left(\int_1^{T^t }\frac{dr}{X^\kappa(r)}- \frac{t\log T}{\kappa p}\right); t\geq 0\right) \Rightarrow (vW(t); t \geq 0)\,. 
\end{align}
The entrance law is given in Remark 4.9 (2) in \cite{patie1}. Let us notice that the case $\delta = 1$ corresponds to a critical continuous state branching process conditioned never to be extinct as mentioned in Remark 4.9 (1) in \cite{patie1}.

Now, to get an invariance principle under $\mathbb Q_a$ for $a> 0$, we have a problem since
 we cannot choose $m$ such that $f_m$ satisfies (\ref{pr}). Nevertheless 
 there is another way to get an invariance principle under $\mathbb Q_a$ for $a> 0$.

We introduce the OU process defined by
\begin{align}
\label{nOU}
\widetilde U(t) := e^{-\kappa^{-1} t} X\left(e^t-1\right)
\end{align}
which is  
 a CBI with immigration mechanism $\chi$ and branching mechanism 
\[\widetilde\phi(\lambda) = \phi(\lambda) + \kappa^{-1} \lambda\,,\]
(see \cite{patie1} section 5.1). 
Let us stress that it is not stationary. We observe that
\[\int_1^\infty \log z \ \mu(dz) < \infty\,,\]
so that applying \cite{friesen2019exponential} Th. 7.7 and Cor. 5.10, we conclude that $\tilde U$ is exponentially ergodic and the  convergence (\ref{C1}) holds under the conditions (\ref{-2}).

Pushing forward this result to the process $X$ we obtain the following Proposition which is an invariance principle for the clock of CBI.
\begin{proposition}
For $a >0$, under $\mathbb Q_a$ as  
$T\rightarrow \infty$,
\begin{align}
\left((\log T)^{-1/2} \left(\int_0^{T^t - 1}\frac{dr}{X^\kappa(r)}- \frac{t\log T}{\kappa p}\right); t\geq 0\right) \Rightarrow (vW(t); t \geq 0)\,. 
\end{align}
\end{proposition}

\begin{remark}
 In \cite{friesen2019exponential} Cor.  5.10, the authors mentioned that  if one starts from a general test function,  it is unlikely to find an explicit formula for the  asymptotic variance in terms of its admissible parameters, except when $f(x) = \exp (\lambda x)$. Our result provides one more example, $f(x) =x^{-\kappa} - (\kappa p)^{-1} $ of a test function with explicit asymptotic variance.  
 \end{remark}

{\it Aknowledgment. \ } This paper was written during a stay of the second author to UNAM during December 2019. He thanks the probability team for its warm hospitality.

\bibliographystyle{plain}
\bibliography{abPM}
\end{document}